\documentstyle[12pt,twoside]{amsart}

\title
{The Nash problem   on  arc families of singularities}
\author{Shihoko Ishii} 
\author{  J\'anos Koll\'ar} 
\address{Shihoko Ishii: Department of Mathematics, Tokyo Institute of
Technology, Oh-Okayama, Meguro, Tokyo, Japan
\newline
e-mail : shihoko@@math.titech.ac.jp}
\address{J\'anos Koll\'ar: Princeton University, Princeton NJ 08544-1000 
USA
\newline
e-mail : 
 kollar@@math.princeton.edu}

\newcommand{\bZ}{{\Bbb Z}}

\newcommand{\bR}{{\Bbb R}}
\newcommand{\bN}{{\Bbb N}}

\newcommand{\Spec}{\operatorname{Spec}}

\newcommand{\Hom}{\operatorname{Hom}}
\newcommand{\ord}{\operatorname{ord}}
\newcommand{\vol}{\operatorname{vol}}
\newcommand{\st}{{\Spec k[[t]]}}
\newcommand{\sT}{{\Spec K[[t]]}}
\newcommand{\D}{{\Delta}}
\let \cedilla =\c
\renewcommand{\c}[0]{{\mathbb C}}  

\renewcommand{\o}[0]{{\mathcal O}} 

\renewcommand{\a}[0]{{\mathbb A}}

\newcommand{\p}[0]{{\mathbb P}}

\newcommand{\q}[0]{{\mathbb Q}}
\newcommand{\map}[0]{\dasharrow}
\newcommand{\qtq}[1]{\quad\mbox{#1}\quad}
\newcommand{\spec}[0]{\operatorname{Spec}}
\newcommand{\sing}[0]{\operatorname{Sing}}
\newcommand{\fs}{{\Spec k[[s]]}}
\newcommand{\fst}{{\Spec k[[s,t]]}}

\def\into{\DOTSB\lhook\joinrel\rightarrow}
\def\to {\longrightarrow}

\newtheorem{thm}{Theorem}[section]

\newtheorem{lem}[thm]{Lemma}
\newtheorem{cor}[thm]{Corollary}
\newtheorem{prop}[thm]{Proposition}
\newtheorem{problem}[thm]{Problem}

\theoremstyle{definition}
\newtheorem{defn}[thm]{Definition}

\newtheorem{say}[thm]{}
\newtheorem{exmp}[thm]{Example}

\newtheorem{rem}[thm]{Remark}

\theoremstyle{remark}

\begin{document}

\begin{abstract}
Nash proved that  every  irreducible component of the
space of arcs through a singularity corresponds to an exceptional
divisor that appears on every resolution. He asked if the converse also 
holds: does  every
such exceptional divisor correspond to an arc family?
We prove that the converse
 holds for toric singularities but fails in general.
\end{abstract}

\maketitle
\markboth{\hfill SHIHOKO ISHII \& J\'ANOS KOLL\'AR\hfill}{\hfill
 THE NASH PROBLEM  \hfill}
 

\section{Introduction}
\noindent
In a 1968 preprint, later  published as 
\cite{nash}, 
 Nash introduced arc spaces and jet schemes for algebraic and
analytic varieties. The problems raised by Nash were studied by 
 Bouvier,  Gonzalez-Sprinberg, Hickel,
Lejeune-Jalabert, Nobile,  Reguera-Lopez and others,
see \cite{B, GL, H, L80, L, LR,  nobile, RL}.
  
  The study of these spaces was further developed by 
   Kontsevich,  Denef and  Loeser 
 as the theory 
  of motivic integration, see \cite{ko, DL}. 
Further interesting applications of jet spaces 
 are given by Musta\cedilla{t}\v{a} \cite{must}.

  The main subject of the paper of Nash is the map from the set of
irreducible components
  of the space of arcs through singular points 
(families of arcs in the original terminology of Nash)
to the set of essential 
  components of a resolution of  singularities.
Roughly speaking, these are the irreducible components of the exceptional 
set  of a given resolution that appear on every
possible resolution,
see Definition \ref{nashessential}.

We call this map the {\it Nash map}, see Theorem \ref{n} for a precise
definition.
  The Nash   map is always injective 
 and Nash asked if it  is always bijective.
This problem remained open  even for  2-dimensional
  singularities, though many cases were settled in \cite{RL}.

  In this paper we prove that the Nash map
is bijective for  toric singularities in  any  dimension,
see Theorem \ref{toric.main.thm}.
  On the other hand we also  show  that the Nash map
is not bijective in general. For instance,
the 4-dimensional hypersurface singularity
$$
x_1^3+x_2^3+x_3^3+x_4^3+x_5^6=0
$$
has only 1 irreducible family of arcs but 2 essential  divisors
over any algebraically closed  field of characteristic $\neq 2,3$.
See Example \ref{bestexmp}.
  
  In \S 2 we define the Nash map and show its injectivity.
This is essentially taken from \cite{nash} with some  scheme
  theoretic details filled in.
  The Nash map for  toric 
  singularities is studied in \S 3.
Counter examples are given in 
   \S 4.  
  
  In this paper, the ground field \( k \) is an algebraically closed 
  field of arbitrary characteristic. A \( k \)-scheme is not 
  necessarily of
  finite type unless we state 
  otherwise.
   A variety means a separated, irreducible and
reduced scheme of finite type over \( 
   k \).
  Every  variety $X$ that we consider is assumed to have a resolution
of singularities $f:Y\to X$ which is an isomorphism
over the smooth locus and whose exceptional set is purely
one codimensional.
Without this or similar  assumptions the definition of essential 
divisors and components
would not make sense.
The existence of resolutions is known in characteristic zero and
for toric varieties in any characteristic.

  The first author would like to thank Professor G\'erard Gonzalez-Sprinberg
  who generated her interest in this problem, provided the information
  on his joint paper \cite{G-S}
 and gave constructive comments to improve this paper.   
We thank   the referee  for useful comments and corrections.
    Part of the work  was completed during the second author's  stay
at the Isaac Newton Institute for Mathematical Sciences, Cambridge.
The first author is partially supported by Grant-in-Aid for
Scientific Research, Japan.
Partial financial support for the second author
was provided by  the NSF under grant numbers 
DMS-9970855 and DMS02-00883.

\section{The space of arcs and the Nash problem}

\begin{defn} Let $X$ be a variety, \( g:X_1 \to X \)
a proper birational morphism from a normal variety \( X_{1} \)
 and $E\subset X_1$ 
an irreducible exceptional divisor  of 
  \( g \). 
 Let  \( f:X_2 \to X \) be another
 proper birational morphism from a normal variety \( X_{2} \).
The birational  map
  \(f^{-1}\circ g  : X_1  \dasharrow  X_2 \) is 
defined on a (nonempty) open subset $E^0$ of $E$.
The closure of $(f^{-1}\circ g)(E^0)$ is well defined.
It is called the {\it center} of $E$ on $X_2$.

  We say that \( E \) appears in \( f \) (or in \( X_2 \)),
  if 
the center of $E$ on $X_2$ is also a divisor. In this case
the birational  map
  \(f^{-1}\circ g  : X_1  \dasharrow  X_2 \) is  a local isomorphism at the 
  generic point of \( E \) and  we denote the birational
  transform of \( E \) on \( X_2 \)
  again by \( E \). For our purposes $E\subset X_1$ is identified
with $E\subset X_2$. (Strictly speaking, we should be talking about the
corresponding {\it divisorial valuation} instead.)
 Such an equivalence class  is called an {\it exceptional divisor over $X$}. 
\end{defn}

\begin{defn}
  Let \( X \) be a variety over \( k \). In this paper, by  a
{\it  resolution} of the 
  singularities of \( X \) 
   we mean  a proper, birational  morphism \( f:Y\to X \) 
with  \( Y \)  non-singular
such that $Y\setminus f^{-1}(\sing X)\to X\setminus \sing X$
 is an isomorphism.

  A resolution  \( f:Y\to X \) is called a {\it divisorial resolution} 
 of \( X \) 
  if 
 the exceptional set is of pure codimension one. 

If $X$ is factorial (or at least $\q$-factorial)
then every resolution is divisorial.
\end{defn}

\begin{defn}\label{nashessential}
  An exceptional divisor \( E \) over \( X \) is called an {\it essential 
  divisor} over \( X \)
   if for every resolution \( f:Y\to X \)
the center of \( E \) on \( Y \) is an irreducible component of 
   \( f^{-1}(\sing X) \).

An exceptional divisor \( E \) over \( X \) is called a {\it 
divisorially essential 
  divisor} over \( X \)
  if  for every divisorial resolution \( f:Y\to X \)
the center of \( E \) on \( Y \) is a divisor,
(and hence also an irreducible component of 
   \( f^{-1}(\sing X) \)).
  
  For a given resolution \( f:Y\to X \),
  the set 
$$
{\cal E}={\cal E}_{Y/X}=
\left\{ 
    \begin{array}{c}
    \mbox{irreducible components of \( f^{-1}(\sing X)  \)}\\  
    \mbox{ which are  centers  of  essential divisors  over \( X \)}
    \end{array} 
    \right\}
$$
corresponds 
  bijectively to
  the set of all essential divisors  over \( X \).

  Therefore we call an element of  \( {\cal E} \) an {\it essential 
  component} on \( Y \). 

Similarly, we can talk about   {\it divisorially essential 
  components} on \( Y \).

  It is clear that  an essential divisor is also a divisorially essential 
  divisor.  We do not know any examples when the two notions are
different. In \S 3 we will see that they coincide  for  
  toric singularities.
\end{defn}

\begin{exmp}
  Let \( (X,x) \) be a normal 2-dimensional singularity.
  Then the set of  the divisorially essential divisors
 over \( X \) coincides with
  the set of the exceptional curves appearing 
  on the minimal resolution $X'\to X$. These are also the essential components
on $X'$.
\end{exmp}


\begin{exmp}\label{nonruled->essntial} Proposition 4 of 
  \cite{abh}
asserts that if $E$ is an exceptional divisor of a  birational morphism
$Y\to Y'$ with $Y'$ smooth then $E$ is ruled,
that is, $E$ is birational to $F\times \p^1$ for some variety $F$.
As noted by Nash, this implies that any nonruled exceptional divisor
of a resolution $Y\to X$ is essential.
\end{exmp}

\begin{exmp}
  Let \( (X,x) \) be a canonical singularity which admits a crepant
  divisorial  
  resolution. 
  A quite large group of such singularities is known (see, for example, 
  \cite{dais} and the references there).
  Then the set of the essential divisors  over \( X \)
  and also the set of the divisorially essential divisors over \( X \)
   coincide  with
  the set of the crepant exceptional divisors.
  Indeed, a divisorial essential divisor should be one of the crepant exceptional
  divisors because of the existence of a divisorial crepant resolution.
  On the other hand, a crepant exceptional divisor  cannot be contracted
  on a non-singular model of \( X \), because if it could be 
  contracted, the 
  discrepancy of the crepant component would have to be positive.
  This shows that every  crepant divisor is an essential component. 
\end{exmp}

\noindent
\begin{defn}
  Let \( X \) be a scheme of finite type over \( k \)
and $K\supset k$ a field extension.
  A  morphism \( \Spec K[t]/(t^{m+1})\to X \) is called an \( m \)-jet
  of \( X \) and \( \Spec K[[t]]\to X \) is called an {\it arc} of \( X \).
  We denote the closed point of \( \Spec K[[t]] \) by \( 0 \) and
  the generic point by \( \eta \).
\end{defn}

\begin{say}
\label{field}
  Let \( X \) be a scheme of finite type over \( k \).
  Let \( {\cal S}ch/k \) be the category of \( k \)-schemes  
   and \( {\cal S}et \) the category of sets.
  Define a contravariant functor  \( F_{m}: {\cal S}ch/k \to {\cal S}et \)
  by 
$$
 F_{m}(Y)=\Hom _{k}(Y\times_{\Spec k}\Spec k[t]/(t^{m+1}), X).
$$
  Then, \( F_{m} \) is representable by a scheme \( X_{m} \) of finite
  type over \( k \), that is
$$
 \Hom _{k}(Y, X_{m})\simeq\Hom _{k}(Y\times_{\Spec k}
\Spec k[t]/(t^{m+1}), X).
$$ 
   This \( X_{m} \) is called the {\it scheme of \( m \)-jets} of \( X \).
  The canonical surjection \( k[t]/(t^{m+1})\to k[t]/(t^{m}) \)
  induces a morphism \( \phi_{m}:X_{m}\to X_{m-1} \).
  Define \( \pi_{m}=\phi_{1}\circ \cdots \circ \phi_{m}:
  X_{m}\to X \).
  A point \( x \in X_{m} \)  gives an \( m \)-jet \( \alpha_{x}:
  \Spec K[t]/(t^{m+1})
  \to X \) and \( \pi_{m}(x)=\alpha_{x}(0) \),
  where \( K \) is the residue field at \( x \). 
  
  Let \( X_{\infty}=\varprojlim _{m}X_{m} \) and call it the 
{\it space of 
  arcs} of \( X \).
\( X_{\infty} \) is not of finite type over $k$
but it is a scheme, see \cite{DL}.
  Denote the canonical projection \( X_{\infty }\to X_{m} \) by \( 
  \eta_{m} \) and the composite \( \pi_{m}\circ \eta_{m} \) by \( \pi \).
    A point \( x \in X_{\infty} \)  
   gives an arc \( \alpha_{x}:\sT
  \to X \) and \( \pi(x)=\alpha_{x}(0) \), 
  where \( K \) is the residue field at \( x \).

  Using  the representability of \( F_{m} \) 
we obtain the following universal property of $X_{\infty}$:
\end{say} 

\begin{prop}
\label{ft}
  Let \( X \) be a scheme of finite type over \( k \).
  Then
  \[ \Hom _{k}(Y, X_{\infty})\simeq\Hom _{k}(Y\widehat\times_{\Spec k}\st, X) \]
   for an arbitrary \( k \)-scheme \( Y \),
   where \( Y\widehat\times_{\Spec k}\st \) means the formal completion 
   of \( Y\times_{\Spec k}\st \) along the subscheme 
   \( Y\times _{\Spec k} \{0\} \).
\end{prop}

\begin{cor}
\label{subscheme}
  There is a universal family of  arcs
  $$
 X_{\infty}\widehat\times _{\Spec k}\st \to X.
$$
\end{cor}

\begin{defn}\label{good.comp}
 Let $X$ be a $k$-variety with singular locus
   \( \sing X\subset X \).
  Every point \( x \) of  the inverse image 
\( \pi^{-1}(\sing X)\subset X_{\infty} \) 
  corresponds to an  arc 
  \( \alpha_{x}:\sT \to X \) such that \( \alpha_{x}(0)\in \sing X \),
  where \( K \) is the residue field at \( x \).
$\pi^{-1}(\sing X)$ is the space of arcs through $\sing X$.

  Decompose  \( \pi^{-1}(\sing X)  \) into its irreducible
  components  
$$
 \pi^{-1}(\sing X)=(\bigcup_{i \in 
  I}C_{i})\cup(\bigcup_{j\in J} 
  C'_{j}), 
$$
 where the \( C_{i} \)'s are  the components with a point \( 
  x \) corresponding to an arc \( \alpha_{x} \)
  such that \( \alpha_{x}(\eta)\not\in \sing X \), while the \( C'_{j} \)'s are
  the components without such points. We call the \( C_{i} \)'s
the {\it good components} of the space of arcs through $\sing X$.

 The notion of ``arc families'' in \cite{nash} is 
the same as the above concept of  good 
  components.
\end{defn}

The next lemma shows that in characteristic zero
every irreducible component
of \( \pi^{-1}(\sing X)\subset X_{\infty} \)  is good.
This can be viewed as a strong form of
Kolchin's irreducibility theorem \cite[Chap.IV,Prop.10]{kolchin}.
See also \cite{gillet}. It is also interesting to compare this with
the results of \cite{must} according to which
the jet spaces $X_m$ are usually reducible.

\begin{lem} Let $k$ be a field of characteristic zero and
 $X$  a $k$-variety. Then 
 every arc through $\sing X$ is a
specialization  of an arc through $\sing X$
whose generic point maps into $X\setminus \sing X$.
\end{lem}

Proof.
We may assume that $X$ is affine. Pick any arc
$\phi:\spec k'[[s]] \to X$ such that $\phi(0)\in \sing X$.
Let $Y$ be the Zariski closure of  the image of $\phi$.
Then $\o_Y$ is an integral domain and $\phi$
corresponds to   an injection $\o_Y\to k'[[s]]$,
where we can take $k'$ to be algebraically closed.
We are done if $Y\not\subset \sing X$.
Otherwise we write $\phi$ as a specialization in two steps.

 First we prove that 
$\phi$ is a specialization
of an arc  $\Phi: \spec K[[s]] \to Y\subset X$
such that $\Phi(0)$ is the generic point of $Y$.

We have an embedding $k'[[s]]\into k'[[S,T]]$
which sends $s$ to $S+T$. 
It is easy to check that the composite
$$
k'[[s]]\into k'[[S,T]]\to k'[[S,T]]/(T)\cong k'[[S]]
$$
is an isomorphism.
Thus we obtain $\Phi$
as the composite
$$
\o_Y\stackrel{\phi}{\to} k'[[s]]\into k'[[S,T]]\into k'((T))[[S]].
$$
Set $K=k'((T^{1/n}:n=1,2,\dots))$, the algebraic closure of $k'((T))$. 
The closed point of $\spec K[[S]]$ maps to the ideal $\Phi^{-1}(S)$,
but the pull back of $(S)$ to $k'[[s]]$ is already the zero ideal.
Thus the closed  point of  $\spec K[[S]]$ maps to the generic point of $Y$.

Repeatedly cutting with hypersurfaces containing $Y$ we obtain a
subvariety $Y\subset Z\subset X$ such that $\dim Z=\dim Y+1$
and $X$ is smooth along the generic points of $Z$.
Let $n:\bar Z\to Z$ be the normalization and $\bar Y\subset \bar Z$
the preimage of $Y$ with reduced scheme structure.
$\bar Y\to Y$ is finite, surjective, and so
generically  \'etale in characteristic zero.
Thus the arc $\Phi: \spec K[[S]] \to Y$ can be lifted to
$\bar{\Phi}: \spec K[[S]] \to \bar Y$. $\bar Z$ is normal,
so smooth along the generic point of $\bar Y$. Thus
$\bar{\Phi}$ is the specialization of an arc
through $\bar Y$ whose generic point maps to the generic point of $\bar Z$.
Projecting  to $Z$ we obtain $\Phi$ and hence $\phi$
as the specialization  of an arc through $\sing X$
whose generic point maps into $X\setminus \sing X$.\qed

\begin{exmp} Let $k$ have characteristic $p$ and consider
the surface $S=(x^p=y^pz)\subset \a^3$
with singular locus $Y=(x=y=0)$.  The normalization
is $\bar S\cong \a^2$ with $(u,v)\mapsto (uv,v,u^p)$. The preimage of $Y$ is
 $\bar Y=(v=0)$  and $\bar Y\to Y$ is purely inseparable.
Thus a smooth arc in $Y$ can not be lifted to $\bar Y$
and it is also not 
the specialization  of an arc through $Y$
whose generic point maps into $S\setminus \sing S$.
In this case  \( \pi^{-1}(\sing S)\subset S_{\infty} \)
has a  component which is not good.
\end{exmp}

\begin{lem}
\label{lift}
  Let \( f:Y\to X \) be a   resolution 
  of the singularities of $X$
and \( E_{1},\ldots,E_{r} \)
  the irreducible components of exceptional sets on \( Y \).
For a good component  \( C_{i} \), let \( C_{i}^o \)
  denote the open subset of \( C_{i} \) consisting of  
 arcs \( \alpha_{x}:\Spec 
   K[[t]]\to X \) such that \( \alpha_{x}(\eta)\not\in \sing X \).
  Then, for every \( x\in C_{i}^o \) the  arc
  \( \alpha_{x} \) can be uniquely lifted to an arc 
  \( \tilde \alpha_{x}:\sT \to Y \).
\end{lem}

\begin{pf}
   As \( f \) is isomorphic outside of \( \sing X \) and \( 
   \alpha_{x}(\eta)\not\in \sing X \),
   we obtain the commutative diagram
   \[ \begin{array}{ccc}
   \Spec K((t))& \to & Y\\
   \downarrow & & \  \downarrow f\\
   \sT & \stackrel{\alpha_{x}}\longrightarrow & X.\\
   \end{array} \]
  Since \( f \) is proper, there exists a unique morphism 
  \( \tilde \alpha_{x}:\sT\to Y \) such that \( f\circ \tilde 
  \alpha_{x}=\alpha_{x} \) by the valuative criterion of properness.   
\end{pf}

  This \( \tilde\alpha_{x} \) is called the lifting of \( \alpha_{x} \).
  Now we have a map 
$$
 \varphi: \mbox{points of }(\bigcup_{i}C_{i}^o)
\to 
\mbox{points of }(\bigcup_{l}E_{l})
$$
given by  
  \( x \mapsto \tilde \alpha_{x}(0) \).
  We emphasize that this map is not  a continuous map of schemes.
  In fact, the image of an irreducible subset is not necessarily 
  irreducible.

\begin{thm}[Nash \cite{nash}]
\label{n}
Let $X$ be a $k$-variety 
and $f:Y\to X$ a resolution of singularities.
Let  \(\{ C_{i}:  i\in I\} \) be 
the good components of the space of arcs through $\sing X$
and  let \( z_{i} \) denote  the generic point of \( C_{i} \).
  Then:
  \begin{enumerate}
  \item[(i)]
  \( \varphi(z_{i}) \) is the generic point of an exceptional component  \( 
  E_{l_{i}}\subset Y \) for some \( l_{i} \).
  \item[(ii)]
  For every \( i\in I \), \( E_{l_{i}} \) is an
 essential component on \( Y \).
  \item[(iii)]
    The resulting  Nash map
$$
\left\{
\begin{array}{c}
\mbox{good components}\\
\mbox{of the space of arcs}\\
\mbox{through $\sing X$}
\end{array}
\right\}
 \stackrel{\cal N}\longrightarrow   
    \left\{ 
    \begin{array}{c}
    \mbox{essential}\\  
    \mbox{components}\\
    \mbox{ on $Y$}
    \end{array} 
    \right\}
\simeq
  \left\{ 
    \begin{array}{c}
    \mbox{essential}\\  
    \mbox{divisors}\\
    \mbox{ over $X$}
    \end{array} 
    \right\}
$$

\noindent
given by   \( C_{i}\mapsto E_{l_{i}} \) is injective.
  In particular, there are only finitely many 
good components of the space of arcs through $\sing X$.
\end{enumerate}
\end{thm}
  
\begin{pf}
  The  resolution \( f:Y\to X  \) induces a morphism
  \( f_{\infty}:Y_{\infty}\to X_{\infty}  \) of schemes. 
  Let \( \pi^Y:Y_{\infty} \to Y \) be the canonical projection.
  As \( Y \) is non-singular, \(  (\pi^Y)^{-1}(E_{l})  \) is 
  irreducible for every \( l \).  
  Denote  by \(  (\pi^Y)^{-1}(E_{l})^o  \) 
  the open subset of \( (\pi^Y)^{-1}(E_{l}) \) consisting of the
   points \( y \) corresponding to   arcs \( \beta_{y}:\Spec 
   K[[t]]\to Y \) such that \( \beta_{y}(\eta)\not\in \bigcup_{l}E_{l} \).
  By  restriction \( f_{\infty} \) gives  \( f'_{\infty}:
  \bigcup_{l=1}^r (\pi^Y)^{-1}(E_{l})^o \to \bigcup_{i\in I}C_{i}^o\).
  For a point \( x\in C_{i}^o \), 
  let \( \alpha_{x}:\sT\to X \) be the corresponding arc,
  where \( K \) is the residue field at \( x \).
  The lifting \( \tilde \alpha_{x}:\sT\to Y \) of \( \alpha_{x} \)
  obtained in Lemma \ref{lift} determines a \( K \)-valued point
  \( \beta:\Spec K\to Y_{\infty} \).
  Denote the image of \( \beta \) by \( y \),
  then \( f_{\infty}(y)=x \).
  Therefore, \( f'_{\infty} \) is surjective.
  Hence,        
 for each \( i\in I 
  \) there is \( 1\leq l_{i}\leq r \) such that the generic point \( 
  y_{l_{i}} \) of
  \( (\pi^Y)^{-1}(E_{l_{i}})^o  \) is mapped to the generic point \( z_{i} \)
  of \( C_{i}^o \).
  Let \( \tilde \alpha_{z_{i}} \) be the lifting of the arc \( 
  \alpha_{z_{i}} \) corresponding to \( z_{i} \) and let 
  \( \beta_{y_{l_{i}}} \) be
  the arc of \( Y \) corresponding to \( y_{l_{i}} \).
   Then \( \beta_{y_{l_{i}}}=\tilde \alpha_{z_{i}} \).
   This is proved as follows:  
  Let \( L \) and \( K \) be the residue fields at \( y_{l_{i}} \) 
  and \( z_{i} \), respectively and \( g:\Spec L[[t]]\to \sT \) be 
  the canonical morphism induced from the inclusion \( K\to L \).
  Then \( \beta_{y_{l_{i}}}=\tilde \alpha_{z_{i}}\circ g \).
  From this, we have \( K=L \) and therefore 
  \( \beta_{y_{l_{i}}}=\tilde \alpha_{z_{i}} \).  
  Note that \( \beta_{y_{l_{i}}}(0)=\pi^Y(y_{l_{i}}) \), which is the generic 
  point  of \( E_{l_{i}} \).
  To finish the proof of (i), just recall 
  \( \varphi(z_{i})=\tilde \alpha_{z_{i}}(0)=
   \beta_{y_{l_{i}}}(0)\).

  Next, we can see that the map  \( C_{i}\mapsto E_{l_{i}} \)
  is injective. 
  Indeed, if \( E_{l_{i}}=E_{l_{j}} \) for \( i\neq j \),
  then \( z_{i}=f'_{\infty}(y_{l_{i}})= f'_{\infty}(y_{l_{j}})= z_{j}\),
  a  contradiction.
  
  To prove that the \( \{E_{l_{i}}: i\in I\} \) are essential 
  components on \( Y \), let \( Y'\to X \) be another  resolution and
   \( \tilde Y \to X \) a divisorial resolution which
  factors through both \( Y \) and \( Y' \).
  Let \( E'_{l_{i}}\subset Y' \) and \( \tilde E_{l_{i}}\subset \tilde Y  \)
   be the irreducible components of the exceptional sets  
  corresponding to \( C_{i} \).
  Then, we can see that \(  E_{l_{i}}  \) and \(  E'_{l_{i}} \) are 
  the image of \( \tilde E_{l_{i}} \).
  This shows that \( \tilde  E_{l_{i}} \) is an essential divisor 
  over \( X \) and therefore \(  E_{l_{i}} \) is an essential 
  component on \( Y \).  
\end{pf}

  Nash poses the following problem in his paper \cite[p.36]{nash}.
  
\begin{problem}
  Is    the Nash map  bijective?
\end{problem} 


\section{The Nash problem for  toric singularities}

\begin{say}
\label{notation}
  We use the notation and terminology of \cite{fulton}.
  Let  $M$  be the free abelian group  ${\bZ}^n$ $(n\geq 2)$
  and  $N$   its dual $\Hom_{\bZ}(M, {\bZ})$.
  We denote  $M\otimes _{\bZ}{\bR}$  and $N\otimes_{\bZ}{\bR}$  by
  ${M_{\bR}}$  and  $N_{\bR}$, respectively.
  The canonical pairing \( \langle\ , \ \rangle:
  N\times M \to \bZ \)  extends to  
  \( \langle\ , \ \rangle:
  N_{\bR}\times M_{\bR} \to \bR \).
  For a finite fan ${\D}$ in ${N_{\bR}}$, the corresponding 
   toric variety is denoted by  $X=X({\D})$.
  For the primitive vector \( v \) in a one-dimensional cone
  \( \tau\in \D \), denote the invariant divisor 
  \( \overline{orb(\tau)} \) in \( X \) by \( D_{v} \).
 
  For a cone \( \tau\in \D \) 
  denote by $U_{\tau}$  the invariant affine open subset which contains
  ${orb\ 
  {\tau}}$ as the unique closed orbit.
   A cone \( \tau \) is called regular or non-singular,
  if its generators can be extended to a basis of \( N \).
  A cone is called singular, if it is not regular.  
  Note that a cone \( \tau \) is regular, if and only if \( U_{\tau} \) is non-singular.
   A cone generated by \( v_{1},\ldots, v_{r} \in N \)
   is denoted by \( \langle  v_{1},\ldots, v_{r} \rangle \).

   We can write $k[M]$  as  $k[x^{u}]_{u\in M}$,
  where we use the shorthand $x^{u}=x_1^{u_1}x_2^{u_2}\cdots x_n^{u_n}$
  for ${u}=(u_1,\ldots ,u_n)\in M$.
\end{say}

\begin{defn}
  An exceptional divisor \( E \) over a toric variety \( X \)
  is called a {\it toric divisorially essential divisor} over \( X \) 
  if its center is a divisor on every
  equivariant divisorial resolution of the singularities of \( X \). 
\end{defn}

  The following is obvious by the definition.
  
\begin{prop}
  For a toric variety \( X \) a divisorially essential divisor over \( X \)
  is a toric divisorially essential divisor over \( X \).
\end{prop}

At this moment the converse of the above proposition is not clear.
But later on, as a corollary of our theorem we obtain the converse.
 
\begin{say}  
  In what follows, we consider  an affine toric variety \( X=X({\D}) \),
  therefore  the fan \( \D \) consists of
  all faces of a cone \( \sigma \).
  Let \( \sigma=\langle e_{1},\ldots,e_{s}\rangle \), where the right 
  hand side means the cone generated by  primitive vectors \( 
  e_{1},\ldots,e_{s} \).
  Let \( T \) be the open orbit in \( X \).
  Let \( W \) be the singular locus of \( X \), then \( 
  W=\bigcup_{\tau:{\operatorname{singular}}}orb(\tau) \).
  Let \( S=N\cap (\bigcup_{\tau:\operatorname{singular}}\tau^o) \),
  where \( ^o \) means the relative interior.    
\end{say}

\begin{prop}
  If \( D_{v} \) is a toric divisorially 
essential divisor  for \( v\in N\cap \sigma \),
  then \( v \) belongs to \( S \).
\end{prop}

\begin{pf}
  If \( D_{v} \) is a toric divisorially 
essential divisor, then the image of \( D_{v} \) must be
  in the singular locus \( W= 
  \bigcup_{\tau:{\operatorname{singular}}}orb(\tau) \).
  Therefore \( v\in S \).
\end{pf}
  
\begin{say}[Sketch of the proof]\label{diagram}
  We  prove that all maps in the following diagram are injective
   and that 
  composite of  all maps is the identity.
  This shows that all maps are bijective.
$$
\begin{array}{ccc}
\left\{
\mbox{minimal elements in $S$}
\right\}
&
  \stackrel{\cal F}{\longrightarrow} 
&
\left\{
\begin{array}{c}
\mbox{good components of}\\
\mbox{arcs through $\sing X$}
\end{array}
\right\} \\

&& \hphantom{{\cal N}}\downarrow {\cal N}\\
{\cal G}\uparrow\hphantom{{\cal G}}
&
 
&
\left\{
\begin{array}{c}
\mbox{essential divisors}\\
\mbox{over $X$}
\end{array}
\right\}\\
 & &  \mbox{ $\cap$}\\
\left\{
\begin{array}{c}
\mbox{toric divisorially}\\
\mbox{essential divisors}\\
\mbox{ over $X$}
\end{array}
\right\} 
&
\supset
&
\left\{
\begin{array}{c}
\mbox{ divisorially essential}\\
\mbox{ divisors over $X$}
\end{array}
\right\}\\
\end{array}
$$
\end{say}

  First we define an order in \( N\cap \sigma \).
  
\begin{defn}
  For two elements \( v, v'\in N\cap \sigma \) we define
  \( v\leq v' \), if \( v'\in v+\sigma \).
  
  For a subset \( A\subset N\cap \sigma \), \( a\in A \) is called 
  minimal in \( A \), if there is no other element \( a'\in A \) such 
  that \( a'\leq a \).
\end{defn} 

  Note that \( v\leq v' \) if and only if \( \langle v, u\rangle\leq 
  \langle v', u\rangle \) for every \( u\in M\cap \sigma^{\vee} \).

  It is clear that \( \leq \) is a partial order, i.e.,
  \begin{enumerate}
  \item[(1)]
  \( v\leq v \),
  \item[(2)]
  if \( v\leq v' \) and \( v'\leq v \),
  then \( v=v' \),
  \item[(3)]
  if \( v\leq v' \) and \( v'\leq v''  \), then
  \( v\leq v'' \).
\end{enumerate}

\begin{defn}
  For an arc \( \alpha:\sT \to X \)  such that \( \alpha(\eta)\in T 
  \),
  define \( v_{\alpha}\in N\cap\sigma \) as follows:
 
  By the condition of \( \alpha \),
  we have a commutative diagram of ring homomorphisms:
  \[ \begin{array}{ccc}
  k[M\cap \sigma ^{\vee}]& \stackrel{\alpha^*}\longrightarrow & K[[t]] \\
  \cap& &\cap \\
  k[M]& \stackrel{\alpha^*}\longrightarrow & K((t)).\\
  \end{array}  \]
  The map \( M\to \bZ \), \( u\mapsto \ord (\alpha^*x^u )\) is a 
  group homomorphism, therefore it determines an element \( v_{\alpha}\in N \)
  such that \( \langle v_{\alpha}, u\rangle = \ord (\alpha^*x^u ) \)
  for every \( u\in M \).
  By the commutative diagram it follows that 
  \( v_{\alpha}|_{M\cap \sigma^{\vee}}\geq 0 
  \),
  hence \( v_{\alpha}\in N\cap \sigma  \).
\end{defn}
 
\begin{prop}
\label{key}
\begin{enumerate}
\item[(i)]
  Let \( \alpha \) be an arc of \( X  \) such that \( \alpha(\eta)\in 
  T \) and \( \tau \)  a 
  face of \( \sigma \).
  Then \( \alpha(0)\in orb(\tau) \), if and only if
  \( v_{\alpha}\in \tau^o \).
  In particular, \( \alpha(0)\in T \), if and only if \( v_{\alpha}=0 \).
\item[(ii)]
  Let \( \Sigma \) be a subdivision of the fan \( \D \) and \( f:Y \to X\) be 
  the toric morphism corresponding to this subdivision.
  Then, an arc \( \alpha \) of \( X  \) such that \( \alpha(\eta)\in 
  T \) is lifted to an arc \( \tilde \alpha \) of \( Y \).
  Let \( \tau\in \Sigma \).
  Then, \( \tilde \alpha(0)\in orb(\tau) \), if and only if 
  \( v_{\alpha}=v_{\tilde \alpha}\in \tau^o \). 
\end{enumerate}
\end{prop}

\begin{pf}
  The first statement of (ii) follows immediately from the properness 
  of \( f \) and the condition \( \alpha(\eta)\in T \).
  The second statement of (ii) follows from the result (i) with replacing
   \( X\) by \(U_{\tau} \).

  For the proof of (i) it is sufficient to prove 
  that \( v_{\alpha}\in \tau \) if and only if \( \alpha(0)\in 
  U_{\tau} \),
  because \( \tau^o =\tau \setminus \bigcup_{\tau'}\tau' \) and 
  \( orb (\tau)=U_{\tau}\setminus \bigcup_{\tau'}U_{\tau'} \),
  where the unions are over all the proper faces \( \tau' \) of \( 
  \tau \).
  The condition \( v_{\alpha}\in \tau \) is equivalent to 
  \( \langle v_{\alpha}, u \rangle \geq 0 \) for all \( u\in 
 M \cap  \tau^{\vee} \).
  And this holds if and only if the ring homomorphism 
  \( \alpha^*:k[M\cap \sigma^{\vee}]\to k[[t]] \) can be extended to
  \( k[M\cap \tau^{\vee}]\to k[[t]] \), which is equivalent to that 
  \( \alpha \) factors through \( U_{\tau} \).
  As \( U_{\tau} \) contains \( T \), this is equivalent to that \( 
  \alpha(0)\in U_{\tau} \).
\end{pf}

\begin{prop}
\label{existence}
  For every point \( v\in S \), there exists an arc 
  \( \alpha:\st \to X \) such that \( \alpha(0)\in W \),
  \( \alpha(\eta)\in T \) and \( v=v_{\alpha} \).
\end{prop}

\begin{pf}
  Define the ring homomorphism \( \alpha^*:k[M]\to k((t)) \)
  by \( \alpha^*(x^u)=t^{\langle v, u \rangle} \).
  Then we have the following commutative diagram:
  \[ \begin{array}{ccc}
  k[M\cap \sigma ^{\vee}]& \stackrel{\alpha^*}\longrightarrow & k[[t]] \\
  \cap& &\cap \\
  k[M]& \stackrel{\alpha^*}\longrightarrow & k((t)),\\
  \end{array}  \]
  because \( \langle v, u \rangle\geq 0 \) for every \( u\in M\cap 
  \sigma^{\vee} \).
  Let \( \alpha :\st \to X \) be the morphism corresponding to \( 
  \alpha^* \), 
  then \( v=v_{\alpha} \) and 
   we obtain \( \alpha(\eta)\in T \)  by the diagram.
  On the other hand, as \( v\in S \), there is a singular face \( 
  \tau< \sigma \) such that \( v=v_{\alpha}\in N\cap \tau^0 \).
  By Proposition \ref{key} \( \alpha(0)\in orb(\tau)\subset W \).
\end{pf}

\begin{prop}[Upper semi-continuity]
\label{upper}
  Let \( C \) be a \( k \)-scheme, 
  \( \alpha:C\widehat\times_{\Spec k}\st \to X \) 
  a family of arcs on \( X \) and 
  \( \alpha_{c}: \spec k(c)[[t]] \to X \) 
  the arc induced from \( \alpha \) for each point \( c\in Y \).
  Here \( k(c) \) is the residue field at \( c \).
  Assume \( \alpha_{c}({\eta })\in T \) for every \( c\in C \).
  Then the map \( C\to N\cap \sigma \), \( c\mapsto v_{\alpha_{c}} \) 
  is upper semi-continuous, i.e., for every \( v\in N\cap \sigma \)
  the subset \(U_{v}:=\{c\in C \mid v_{\alpha_{c}}  \leq v \} \) is
  open in \( C\). 
  In particular, 
    if there is  a point \( z\in C \) such that \( v_{\alpha_{z}} \) is 
  minimal in \( S \),
  then there is a non-empty open  subset \( U\subset C \) such that 
  \( v_{\alpha_{c}}=v_{\alpha_{z}} \) holds for every \( c\in U \).

\end{prop}

\begin{pf}
  It is sufficient to prove the assertion in the  affine case \( 
  C=\Spec A \).
  Let \( \alpha^*:k[M\cap \sigma^{\vee}] \to A[[t]] \) be the ring 
  homomorphism corresponding to \( \alpha \).
  Let \( \alpha^*(x^u)  \) be \( a_{0}^u+a_{1}^ut+a_{2}^ut^2+\ldots \),
  where \( a_{i}^u\in A  \) for \( i\geq 0 \).
  By the definition of \( U_{v} \), 
  a point \( c\in C \) belongs to \( U_{v} \), 
  if and only if
  \( \langle v_{\alpha_{c}}, u \rangle  \leq \langle v, u \rangle \)
  for every \( u\in M\cap \sigma^{\vee} \).
  This is equivalent to that 
  for every element \( u \) of generating system of \( M\cap\sigma^{\vee} \) 
  there exists \( i\leq \langle 
  v, u \rangle \) such that \( a_{i}^u(c)\neq 0 \).
  Now, we see that \( U_{v} \) is a finite union of the complements of
   zero locus of  
  functions 
  on \( C \).        
\end{pf}

\begin{say}
  Let \( \{C_{i}: i\in I\} \) be the  good components of the space of
arcs through $W$.
  For each component \( C_{i}\subset X_{\infty} \),
  there exists a corresponding family 
  \(\alpha_{i}: C_{i}\widehat\times_{\spec k} \st \to X \)
  of arcs  
  by Corollary \ref{subscheme}.   
\end{say}

\begin{lem}
\label{non-empty}
 Under the above notation, 
 for a minimal element \( v\in S \) there are a good component
 \( C_{i} \) and a non-empty open subset \( U\subset C_{i} \)
 such that \( v_{\alpha_{ic}}=v \) for every \( c\in U \),
 where \( \alpha_{ic}:\spec k(c)[[t]] \to X \) is th arc induced
 from \( \alpha_{i} \).
 
 For a minimal element \( v\in S \), take one of these components 
 \( C_{i} \) and define \( {\cal F}(v):=C_{i} \).
 Then the map 
 \(  \{
 \mbox{minimal elements in}\ S 
 \} \) \( \stackrel{\cal F}\longrightarrow \)
\( \left\{
C_{i}\\
\right\} \)
is injective.
\end{lem}

\begin{pf}
  For a given minimal element \( v\in S \) there is an arc
  \( \alpha:\st \to X \) such that \( \alpha(0)\in W \),
  \( \alpha(\eta)\in T \) and \( v_{\alpha}=v \) by Proposition 
  \ref{existence}.
  Then, by Definition \ref{good.comp}, there exist a good component \( C_{i} \) 
  and its \( k \)-valued point \( z \) such that \( \alpha=\alpha_{iz} \).
  As \( \alpha_{i}(C_{i}\times_{\Spec k} \{0\})\subset W \) and
  \( \alpha_{iz}(\eta)\in T \),
  there exists a non-empty open subset \( V\subset C_{i} \) such that 
   both the
  conditions \( \alpha_{i}(V\times_{\Spec k} \{0\})\subset W \) and 
  \( \alpha_{ic}(\eta)\in T \)  for every \( c\in V \) hold.
  Then, by Proposition \ref{upper}, there exists a non-empty open subset 
  \( U\subset V \) such that \( v_{\alpha_{ic}}=v \).
  
  The second assertion is obvious from the first statement.      
\end{pf}

\begin{lem}
\label{composite}
  Let \(  {\cal N}: \{C_{i}: i\in I\} \to  \{ \mbox{essential
  divisors}\}\),
  \( C_{i}\mapsto E_{l_{i}}   \) be the Nash map in 
  Theorem \ref{n}.
  Then the composite
$$
  { \cal N\circ F}:
\{\mbox{minimal elements in $S$}\} \to
 \{ \mbox{essential divisors}\}
$$
  satisfies \( {\cal N\circ F}(v)=
 D_{v} \).
\end{lem}

\begin{pf}
  By Lemma \ref{non-empty}, the generic point \( z \) of \( 
  {\cal F}(v) \) corresponds to an arc \( \alpha:\sT\to X \) such that
  \( v_{\alpha}=v \).
  Let \( \tilde \alpha \) be the lifting of \( \alpha \) 
  as an arc of a toric  divisorial  resolution 
  \( Y \). 
  By the definition of \( {\cal N} \), 
  \(  {\cal N\circ F}(v) \) is an exceptional divisor containing
  \( \tilde \alpha(0) \) as the generic point.
  By Proposition \ref{key}, the exceptional divisor 
  \(\overline{orb(\tau)}  \) containing \( \tilde \alpha(0) \)
  satisfies \( v=v_{\alpha}=v_{\tilde \alpha}\in \tau^o \).
  Therefore this divisor is \( D_{v} \).     
\end{pf}

 We prove the following by using the idea of the proof of 
 \cite[Th\'eor\`eme 1.10]{G-S}.

\begin{lem} 
\label{g-s}
  Consider the map
 $$
 {\cal G}:\{\mbox{toric divisorially
 essential divisors over $X$}\} \to S 
$$
 given by \( {\cal G}(D_{v})=v \).
  Then, this map is injective and its image  is 
contained in the set of
 minimal elements of $S$.\qed
\end{lem}

\begin{pf}
  The injectivity is clear by the definition of the map.
  For the second assertion it is sufficient to prove that
  if a primitive vector \( v\in S \) is not minimal then \( D_{v} \) is not
  toric divisorially essential.
To do this, 
  we construct  a regular subdivision \( \Sigma \) of
  \( \sigma \) such that the map \( X(\Sigma)\to X \) is a divisorial 
  resolution of \( X \), and in which \( v\bR_{\geq 0} \) does not appear as a 
  one-dimensional cone.  

  If \( v\in S \) is not minimal, then \( v \) can be  written as 
  \( v=n_{1}+n_{2} \), where \( n_{1}\in S \) and \( n_{2}\in N\cap\sigma
  \setminus \{0\} \).
  Then, we can reduce it into two cases: (1) \( n_{1},n_{2}\in S \), 
  (2) \( n_{1}\in S \) and \( n_{2}  \) is in a one-dimensional face 
  of \( \sigma \).
  Indeed, if \( n_{2}\not\in S \), then \( n_{2} \) is in a 
  non-singular face \( \tau \) of \( \sigma \).
  Let \( \tau=\langle e_{1},\ldots,e_{d}\rangle \), then \( n_{2}=
  \sum_{i=1}^d b_{i}e_{i}\) with \( b_{i}\in \bN\cup \{0\} \) 
  \( (i=1,\ldots,d) \).
  We may assume that \( b_{1}\neq 0 \).
  Let \( \gamma \) be the minimal face of \( \sigma \) containing the 
  cone \(\langle n_{1}, \sum_{i=2}^d b_{i}e_{i}\rangle \),
  then, since \( n_{1}\in \gamma \), \( \gamma \) is singular  and
   \( n_{1}+\sum_{i=2}^d b_{i}e_{i}\in \gamma^o\subset S \).
  Here, replace \( n_{1} \) by    \( n_{1}+\sum_{i=2}^d b_{i}e_{i} \)
   and $n_2$  by  $b_1e_1$, then we can reduce to the case (2).
  
  Next, take the minimal regular subdivision of the 2-dimensional cone
  \( \langle n_{1}, n_{2}\rangle \) 
  (\cite[Proposition 1.8]{G-S}) which gives the minimal resolution of
  the 2-dimensional singularity.
  Let \( \langle v_{1},v_{2}\rangle \) be its 2-dimensional cone 
  containing \( v \), then \( v \) is in the relative interior of this
  cone.
  We will construct a regular subdivision of \( \sigma \) which 
  contains \( \langle v_{1},v_{2}\rangle \) as a cone.
  We may assume that \( v_{1}\in S \).    
  First, take the star-shaped subdivision \( \Sigma_{1} \)
   with the center \( v_{1} \).
  Then take the star-shaped subdivision \( \Sigma _{2} \) of \( 
  \Sigma_{1} \) with the center \( v_{2} \) 
  if  \( v_{1},v_{2} \) are in the case (1).
  If \( v_{1},v_{2} \) are in the case (2), let \( \Sigma_{2}=\Sigma_{1} \).
  Here, we note that the exceptional set for the corresponding 
  equivariant morphism is a divisor.
  If \( \Sigma_{2} \) is not simplicial, let \( \gamma \) be a 
  minimal  dimensional cone which is not simplicial.
  Take \( n\in \gamma^o \) and take the star-shaped subdivision
  of \( \Sigma_{2} \) with
  the center \( n \). 
  Then \( \gamma  \) is divided into simplicial cones and 
  the exceptional set for the corresponding 
  equivariant morphism is a divisor.
  Continuing this procedure, we finally obtain a simplicial subdivision
  \( \Sigma_{3} \).
  If \( \Sigma_{3} \) is not regular, take a cone \( \lambda=
  \langle p_{1},\ldots,p_{t}\rangle\in \Sigma_{3} \) with
  the maximal multiplicity.
  The  multiplicity is \( \vol P_{\lambda} \),
  where \( P_{\lambda}=\{\sum_{i=1}^t c_{i}p_{i}\mid 0\leq c_{i}<1\}\).
  Since \( \vol P_{\lambda}>1 \), there is a non-zero element 
  \(n'\in P_{\lambda}\cap 
  N \).
  Take the star-shaped subdivision with the center \( n' \).
  Then again    the exceptional set for the corresponding 
  equivariant morphism is a divisor.
  Continuing this procedure, we finally obtain a regular subdivision
  \( \Sigma_{4} \).
  As we did not change the cone \( \langle v_{1},v_{2}\rangle \) in 
  these procedures, \( \Sigma_{4} \) contains this cone.
  Therefore, the exceptional divisor \( D_{v} \) does not appear in 
  \( X(\Sigma_{4})  \). 
  As all  regular cones are unchanged, the corresponding equivariant
  morphism is a resolution which is isomorphic outside  the singular locus.
  It is clear that the resolution is divisorial, as we saw it in 
  each step of subdivisions.              
\end{pf}

\begin{thm}\label{toric.main.thm}
  Let \( X \) be  an affine toric variety. Then the  Nash map 
  $$
 {\cal N}: \{C_{i}: i\in I\}\to  \{\mbox{essential
  divisors over $X$} \} 
$$
 is bijective.
\end{thm}

\begin{pf}
  In the diagram \ref{diagram}, we obtain that \( {\cal F} \) is 
  injective by Lemma \ref{non-empty}, \( {\cal N} \) is 
  injective by Nash's theorem \ref{n} and \( {\cal G} \) is injective
  by Lemma \ref{g-s}.
  We also have that \( {\cal G}\circ{\cal N}\circ {\cal F} \) is the 
  identity map on \{minimal elements in \( S \)\} by Lemma 
  \ref{composite} and \ref{g-s}.
  Hence, \( {\cal G}, {\cal N}, {\cal F}  \) are all bijective.
\end{pf}

By the proof of the above theorem, the following are obvious.

\begin{cor}
  For a toric variety \( X \), \( E \) is an essential divisor over \( 
  X \), if and only if \( E \) is a toric divisorially
 essential divisor over \( X \).
\end{cor}

The analogous result for essential divisors is proved in \cite{B},
but the definition used there is not quite equivalent to ours.

\begin{cor}
  For a cone \( \sigma \) in \( N \) the number of the minimal elements in
  \( S=N\cap(\bigcup_{\tau:\operatorname{singular}}\tau^o) \) is finite.
  More precisely this number is the number of essential components 
  and also the number of the good components.
\end{cor}

\begin{cor}
  For a general point \( c\in C_{i} \) for 
\( (i\in I) \), the corresponding arc 
  \( \alpha_{ic} \) satisfies \( \alpha_{ic}(\eta)\in T \).
\end{cor}

\begin{exmp}
  Let \( e_{1}=(1,0,0), e_{2}=(0,1,0), e_{3}=(1,1,e)\in N\simeq \bZ^3 \)
  and \( \sigma=\langle e_{1},e_{2},e_{3}\rangle \).
  Then all proper faces of \( \sigma \) are regular and \( \sigma \) 
  itself is not regular, therefore the  affine toric
  variety \( X \) corresponding to \( \sigma \) has an isolated 
  singularity at the closed orbit.
  We can also see that \( S=N\cap\sigma^o \).
  By simple calculations we obtain that the minimal elements in
  \( S \) are  \( (1,1,d) \) \( (1\leq d \leq e-1) \).
  Therefore, by our theorem the number of  \( C_{i} \)'s 
  and the number of the essential components are both \( e-1 \).      
\end{exmp}

\section{Counter examples to the Nash problem}

\noindent
The basic idea of our counter examples to the Nash problem is  the 
following:

Take a singularity $x\in X$ and a partial resolution
$p:Y\to X$ with exceptional divisor $F\subset Y$.
Assume that  $Y$ has a singular point $y\in F$
such that every general arc $g:\fs \to (Y,y)$
is contained in an embedded smooth 
 surface germ $G:\fst\to (Y,y)$.

Assume that there is an essential divisor $E$ 
over $X$ whose center on $Y$ is $y$.
The arcs
on $Y$ that should correspond to $E$
 are all arcs through $y$. If such an arc is contained in
an embedded smooth surface germ $G:\fst\to (Y,y)$,
 then this arc can be moved in $Y$ such that its closed point moves
along the curve $G^{-1}(F)$, hence the arcs through $E$ are all 
limits of arcs through some component of $F$.

This implies that $E$ does not correspond to an irreducible component
of the family of arcs through $x\in X$.
If we can also arrange $E$ to be essential, we have a counter example to
the Nash problem.

\begin{say}
Algebraically, a smooth formal curve through $0\in Y$ is 
equivalent to a surjection
$\phi:\hat{\o}_Y\to k[[s]]$,  where $\hat{\o}_Y$
denotes the completion of $\o_Y$ at the ideal $m_0$ of $0$.
Similarly, a smooth  surface germ is equivalent to a surjection
  $\Phi:\hat{\o}_Y\to  k[[t,s]]$.
The induced maps $m_0/m_0^2\to (s)/(s)^2$
and $m_0/m_0^2\to (s,t)/(s,t)^2$
correspond to a point  and a line in 
the exceptional divisor of the blow up $B_0Y\to Y$.
\end{say}

\begin{lem} \label{embdisc.lem}
Let $0\in Y\subset \a^n$ be a hypersurface singularity
of multiplicity $m$ defined by an equation $F=0$ where
$F=F_m+F_{m+1}+\dots$ is the decomposition into  homogeneous pieces.
Set $Z=(F_m=0)\subset \p^{n-1}$ and let $z\in Z$
be a point and $z\in L\subset Z$ a line such that $Z$ is smooth along $L$
and $H^1(L,N_{L|Z})=0$.

Let $\phi:\hat{\o}_Y\to k[[s]]$ be a smooth formal curve through $0$
with tangent direction $z$. Then $\phi$ can be extended
to  a surjection $\Phi:\hat{\o}_Y\to  k[[t,s]]$ with tangent direction $L$.
\end{lem}

\begin{pf}
 The line $L$ can be identified with a map
$\Phi_1:k[y_1,\dots,y_n]\to k[s,t]$ such that the
$\Phi_1(y_i)$ are linear in $s,t$ and 
$\Phi_1(F)\in (s,t)^{m+1}$. 
Our aim is to find inductively maps
$$
\Phi_r:k[y_1,\dots,y_n]\to k[s,t]
\qtq{such that} \Phi_r(F)\in (s,t)^{m+r},
$$
$\Phi_r$ modulo $(t)$ coincides with $\phi$ modulo $(s^{r+1})$
and $\Phi_r$ is congruent to $\Phi_{r+1}$ modulo  $(s,t)^{r+1}$.
If this can be done then the  inverse limit of the maps
$$
k[y_1,\dots,y_n]\stackrel{\Phi_r}\to k[s,t]\to k[s,t]/(s,t)^{r+1}
$$
gives $\Phi: k[[y_1,\dots,y_n]]\to  k[[s,t]]$
such that $\Phi(F)=0$. Thus it descends to
 $\Phi:\hat{\o}_Y\to  k[[s,t]]$

A map $g:k[y_1,\dots,y_n]\to (\mbox{any ring})$ can be identified with
the vector $(g(y_1),\dots,g(y_n))$. Using this convention,
by changing coordinates we may assume that
$\phi=(s,0,\dots,0)$ and $L=(y_3=\cdots=y_n=0)$.
The first condition implies that no power of $y_1$
appears in $F$ and the second means that we can choose
 $\Phi_1=(s,t,0,\dots,0)$. 

Assume that we already have
$\Phi_r$ which we assume to be of the form
$$
\Phi_r=(s,t,tA_{3,r-1}(s,t),\dots, tA_{n,r-1}(s,t))
$$
where the $A_{i,r-1}$ are polynomials of degree $\leq r-1$
without constant terms. The vanishing of the constant
term comes from extending the map $\Phi_1$ and the
divisibility by $t$ comes from the requirement of extending  $\phi$.
We are looking for $\Phi_{r+1}$ of the form
$$
\Phi_{r+1}=
(s,t,tA_{3,r-1}(s,t)+tB_{3,r}(s,t),\dots, tA_{n,r-1}(s,t)+tB_{n,r}(s,t)),
$$
where the $B_{i,r}$ are homogeneous of degree $r$.
Let us compute $\Phi_{r+1}(F)$. Using the Taylor expansion, we get that
\begin{eqnarray*}
\Phi_{r+1}(F)&=&\Phi_{r}(F)+
t\cdot \sum_{i=3}^n\frac{\partial F_m}{\partial y_i}(s,t,0,\dots,0)
\cdot  B_{i,r}(s,t)\\
&& +\mbox{(terms of multiplicity $\geq m+r+1$).}
\end{eqnarray*}
By the inductive assumption,
$$
\Phi_{r}(F)= t\cdot C_{m+r-1}(s,t)+\mbox{(terms of multiplicity $\geq m+r+1$)},
$$
where $C_{m+r-1}$ has degree $m+r-1$.
In order to achieve that $\Phi_{r+1}(F)\in (s,t)^{m+r+1}$,
  we need  to find polynomials $B_{i,r}$
such that
$$
C_{m+r-1}(s,t)=-\sum_{i=3}^n\frac{\partial F_m}{\partial y_i}(s,t,0,\dots,0)
\cdot  B_{i,r}(s,t).\eqno{(*)}
$$
Since we know nothing about $C_{m+r-1}$, we need to guarantee that
the ideal generated by 
the partials ${\partial F_m}/{\partial y_i}(s,t,0,\dots,0)$
contains all homogeneous polynomials of degree $m+r-1$ in $s,t$
for every $r\geq 1$.
The critical case is $r=1$.

The normal bundles of $L$ in $Z$ and in $\p^{n-1}$
are related by an exact sequence
$$
0\to 
N_{L|Z}\to N_{L|\p^{n-1}}\cong \o(1)^{n-2} \stackrel{dF_m}{\longrightarrow}
N_{Z|\p^{n-1}}|_L\cong \o(m)\to 0,
$$
and $dF_{m}$ is  the map $\o(1)^{n-2}\to \o(m)$
given by multiplication by the partials ${\partial F_m}/{\partial y_i}$
for $i=3,\dots,n$. We have assumed that $H^1(L,N_{L|Z})=0$,
thus the induced map
\begin{eqnarray*}
dF_m: \sum_{i=3}^nH^0(L,\o(1))&\to & H^0(L,\o(m)),\qtq{given by}\\
(l_3,\dots,l_n)&\mapsto &
 \sum_{i=3}^nl_i\frac{\partial F_m}{\partial y_i}(s,t,0,\dots,0)
\end{eqnarray*}
is surjective.
Thus the equation (*) always has a solution.
\end{pf}

\begin{thm}\label{nashthm}
 Let $Z\subset \p^{n-1}$ be a smooth hypersurface.
Assume that $Z$ is  not  ruled but through
a general point of $Z$ there is  a line $L$
such that $H^1(L,N_{L|Z})=0$.

Let $0\in X$ be any singularity with a partial resolution
$p:Y\to X$  and $y\in p^{-1}(0)$ a point such that
\begin{enumerate}
\item[(1)]  $y\in Y$ is a hypersurface singularity whose 
 projectivised tangent cone  is isomorphic to $Z$, and
\item[(2)] $p^{-1}(0)\subset Y$ is a Cartier divisor.
\end{enumerate}

Then the blow up $B_yY$ gives an essential exceptional divisor
$Z\cong E\subset B_yY$ over $0\in X$ which
 does not correspond to an
irreducible family of arcs on $X$.
\end{thm}

\begin{pf}
$E$ is an essential divisor by Example \ref{nonruled->essntial}.

In order to prove that $E$ does not correspond to an
irreducible family of arcs on $X$,
consider the family  $W$ of arcs in $B_yY$ through $E$.
These correspond to a subset 
$W_y$ of arcs on $Y$ through $y$ and to a subset 
$W_x$ of 
arcs in $X$ through $x$. We claim that
$W_x$ is not an irreducible component of the family of all arcs
on $X$ through $x$.

In order to see this, it is enough to show that a general arc in
$W_y$ is a limit of arcs in $Y$ through $p^{-1}(0)$ 
but not passing through $y$.

By assumption, the pull back of a general local equation of
 $y$ contains $E$ with multiplicity 1.
A general arc in $W$  is transversal to $E$,
so the general member of $W_y$ is an arc on $Y$ 
wich has multiplicity 1 intersection with a general local equation of
 $y$.
Hence the general member of $W_y$ is a smooth  arc on $Y$
with general
 tangent direction. Therefore, by our assumption on $Z$ 
and  by (\ref{embdisc.lem}), a general arc in $W_y$
is contained in a smooth surface germ.
Thus it is a limit of arcs through $p^{-1}(0)$
 which do not pass through $y$.
Hence $W_x$ is not an irreducible component of the space of
arcs on $X$ through $x$. 
\end{pf}

\begin{rem}\label{nonruled} 
In characteristic 0, a smooth
 hypersurface $Z\subset \p^{n-1}$
is covered by lines if and only if $\deg Z\leq n-2$
(cf. \cite[V.4.6]{rcbook}).
A general line then has $H^1(L,N_{L|Z})=0$ by
\cite[II.3.11]{rcbook}.
 Thus the key
 condition is to check that $Z$ is not birationally ruled.
This can not happen if $n\leq 4$. In higher dimensions
there are two
 known  sets of examples:

(1) $Z\subset \p^4$ is a smooth cubic. Then $Z$ is not birational
to $\p^3$.
This was proved by 
\cite{cl-gr} over $\c$ and 
by \cite{murre} in characteristic $\neq 2$.
 This  implies that $Z$ is not ruled.
Indeed, assume that $Z$ is birational to $S\times \p^1$.
There is a degree 2 map $\p^3\map Z$ (this goes back to M.\ Noether, 
cf.\ \cite[V.5.18.3]{rcbook}),
 so in characteristic $\neq 2$ we get a dominant separable  map
$\p^3\map S\times \p^1\to S$. Thus $S$ is 
 rational
 by Castelnuovo's theorem. Therefore $Z$ is birational to $\p^2\times \p^1$
and so rational, a contradiction.

Every line on a smooth cubic satisfies  $H^1(L,N_{L|Z})=0$ by
\cite[V.4.4.1]{rcbook} in any characteristic.

(2) $Z\subset \p^{n-1}$ is a very general hypersurface
with $n\geq \deg Z\geq \tfrac{2n}{3}+2$. These are nonruled 
in characteristic zero by
\cite{nonrat}.
\end{rem}

\begin{exmp}\label{bestexmp} The 4-dimensional hypersurface singularity
over an algebraically closed  field of characteristic $\neq 2,3$ 
$$
x_1^3+x_2^3+x_3^3+x_4^3+x_5^6=0
$$
has only 1 irreducible family of arcs but 2 essential exceptional components.
\end{exmp}

\begin{pf}
Apply  Theorem \ref{nashthm}  to  $X=(x_1^3+x_2^3+x_3^3+x_4^3+x_5^6=0)$.
Blowing up the origin produces $Y$. The exceptional divisor $F\subset Y$ 
is Cartier
and $Y$ has a unique singular point which is the cone over
the cubic 3--fold $Z:=(x_1^3+x_2^3+x_3^3+x_4^3+x_5^3=0)$.
$Z$ is not birationally ruled by  (\ref{nonruled}.1).

Blowing up the unique singular point of $Y$ we get a resolution
of $X$ with 2 exceptional divisors.
One is  $E\cong Z$ and the other is $F'$, the birational transform of $F$.

$F'$ is birationally ruled, but it is still essential.
Indeed, the family of arcs on $X$ has to correspond to
some exceptional divisor, and $F'$ is the only possibility.
Thus $F'$ has to be essential. Another way to see this is to note
that $X$ is terminal and $F'$ has minimal discrepancy, namely 1.
\end{pf}


\makeatletter \renewcommand{\@biblabel}[1]{\hfill#1.}\makeatother


\begin{thebibliography}{11}

\bibitem{abh}  S.\ Abhyankar,
{\em On the valuations centered in a local domain.}
Amer. J. Math. 78 (1956) 321--348. 

\bibitem{artin} M.\ Artin, {\em
On the solutions of analytic equations}, 
Invent. Math. 5 1968 277--291.

\bibitem{B} C. Bouvier, {\em  Diviseurs essentiels, composantes essentielles
 des vari\'et\'es
toriques singuli\`eres},
Duke Math. J. 91 (1998) 609--620

\bibitem{G-S} C. Bouvier and G. Gonzalez-Sprinberg, {\em
Syst\`eme g\'en\'erateur minimal, diviseurs essentiels et 
G-d\'esingularisations de vari\'et\'es toriques}, Tohoku Math. J.
{\bf 47}, (1995) 125--149. 



\bibitem{cl-gr} C.H.\ Clemens and  P.\ Griffiths, {\em
The intermediate Jacobian of the cubic threefold},
Ann. of Math.  95 (1972), 281--356




\bibitem{dais} D. I. Dais, C. Haase and G. M. Ziegler, {\em
All toric local complete intersection singularities admit projective crepant
resolutions}, Tohoku Math. J. (2) {\bf 53},  (2001) 95--107.


\bibitem{DL} J. Denef and F. Loeser, {\em Germs of arcs on singular
varieties and motivic integration,} Invent. Math. {\bf 135}, (1999)
201--232. 
 
\bibitem{fulton} W. Fulton, {\em Introduction to Toric Varieties},
Annals of Math. St. {\bf 131}, (1993) Princeton University Press.

\bibitem{gillet} H. Gillet, {\em 
Differential Algebra - a scheme theory approach}, preprint,
http://www.math.uic.edu/~henri


\bibitem{GL} G. Gonzalez-Sprinberg and M. Lejeune-Jalabert, {\em Families of smooth curves on
surface singularities and wedges}, Annales Polonici Mathematici, 
LXVII.2, (1997) 179--190.

\bibitem{H} M. Hickel, {\em Fonction de Artin et germes de courbes 
trac\'ees sur un germe d'espace analytique}, Amer. J. Math. {\bf 115},
(1993) 1299--1334.

\bibitem{kolchin} E. R.  Kolchin
{\em Differential algebra and algebraic groups},
 Pure and Applied Mathematics, Vol. 54. Academic Press, New York-London, 1973.

\bibitem{nonrat} J.\ Koll\'ar, {\em
Nonrational hypersurfaces}, J. Amer. Math. Soc. 8 (1995)
241--249

\bibitem{rcbook} J.\ Koll\'ar, {\em
Rational Curves on Algebraic Varieties}, Springer, 1996


\bibitem{ko} M. Kontsevich, {\em Lecture at Orsay} (December 7, 1995)

\bibitem{L80} M. Lejeune-Jalabert, {\em Arcs analytiques et
r\'esolution minimale des 
surfaces quasihomog\`enes}. in: Lecture Notes in  Math. {\bf 777}, (1980) 
303--336.

\bibitem{L} M. Lejeune-Jalabert, {\em Courbes trac\'ees sur un germe 
d'hypersurface}. Amer. J. Math. {\bf 112}, (1990) 525--568.

\bibitem{LR} M. Lejeune-Jalabert and A. J. Reguera-Lopez, {\em Arcs 
and wedges on sandwiched surface singularities}, Amer. J. Math. 
{\bf 121}, (1999) 1191--1213. 

\bibitem{murre}  J. P. Murre, 
{\em Reduction of the proof of the non-rationality of a
 non-singular cubic threefold to a result of
Mumford},
Compositio Math. 27 (1973), 63--82

\bibitem{must} M. Musta\cedilla{t}\v{a}, {\em
Singularities of Pairs via Jet Schemes},
 J. Amer. Math. Soc. 15 (2002), 599-615.

\bibitem{nash} J. F. Nash, {\em Arc structure of singularities},
Duke Math. J. {\bf 81}, (1995) 31--38. 

\bibitem{nobile} A. 
Nobile, {\em On Nash theory of arc structure of singularities}. Ann. Mat. Pura
Appl.  160 (1991), 129--146

\bibitem{RL} A. J. Reguera-Lopez, {\em Families of arcs 
on rational surface singularities}, Manuscr. Math. 
{\bf 88}, (1995) 321--333. 
\end{thebibliography}
\end{document}